\documentclass[11pt]{amsart}

\usepackage{graphicx}    
\usepackage{amssymb}
\usepackage{amsthm}
\usepackage{amsmath}
\usepackage[T1]{fontenc}
\theoremstyle{definition} \newtheorem{prop}{Proposition}
\theoremstyle{definition} \newtheorem{Lemma}{Lemma}
\theoremstyle{definition} \newtheorem{defi}{Definition}
\theoremstyle{definition} \newtheorem{coro}{Corollary}

\DeclareMathOperator{\LMAX}{\textsc{lmax}}
\DeclareMathOperator{\RMAX}{\textsc{rmax}}
\DeclareMathOperator{\Av}{Av}
\DeclareMathOperator{\Ins}{Ins}
\DeclareMathOperator{\EGF}{\textsc{egf}}

\title[Patterns of the form $k$-$\sigma$-$k$]{Avoidance of Partially Ordered Generalized Patterns of the form $k$-$\sigma$-$k$}
\author[\textrm{Marteinn T. Hardarson}]{\textrm{Marteinn T. Hardarson} \\ \\ \textit{Reykjavík University}}

\begin{document}

\bibliographystyle{abbrv} 

\begin{abstract}
Sergey Kitaev \cite{kita}
has shown that the exponential generating function for
permutations avoiding the generalized pattern $\sigma$-$k$, where $\sigma$ is a pattern without dashes
and $k$ is one greater than the biggest element in $\sigma$, is determined by the exponential generating
function for permutations avoiding $\sigma$.

We show that this also holds for permutations avoiding all the generalized patterns
$\sigma_1$-$k_1$, $\dots$, $\sigma_n$-$k_n$, where $\sigma_1$, $\dots$, $\sigma_n$ are patterns without
dashes and $k_i$ is one greater than the biggest element in $\sigma_i$. Similarly the
exponential generating function for permutations avoiding the partially ordered generalized patterns
$k_1$-$\sigma_1$-$k_1$, $\dots$, $k_n$-$\sigma_n$-$k_n$ can be determined from the exponential
generating function for permutations avoiding the generalized patterns $\sigma_1$, $\dots$, $\sigma_n$,
where $\sigma_1$, $\dots$, $\sigma_n$ are patterns without dashes and $k_i$ is one greater than the
largest element in $\sigma_i$. Since $k$ is the greatest element in the pattern $k$-$\sigma$-$k$,
avoidance of $k$-$\sigma$-$k$ is equivalent to simultaneous 
avoidance of $(k+1)$-$\sigma$-$k$ and $k$-$\sigma$-$(k+1)$.

Using this we construct a bijection between bicolored set partitions and permutations avoiding the
partially ordered generalized pattern $3$-$12$-$3$
(that is, permutations avoiding both the patterns $3$-$12$-$4$ and $4$-$12$-$3$).
By using this method twice, we find a closed formula for the exponential generating function for permutations
avoiding the partially ordered generalized pattern $3$-$121$-$3$.

Finally, we give a complete classification of when single partially ordered generalized patterns have the same 
set of avoiders.

\end{abstract}

\maketitle

\thispagestyle{empty}

\newpage

\section{Introduction}

We write a classical pattern $P$ as a permutation $p_1p_2\dots p_m$.
Let $\pi = a_1a_2\dots a_n$ be a permutation and $1 \le k_1<k_2<\dots<k_m \le n$.
We say that $a_{k_1}a_{k_2}\dots a_{k_m}$ is an occurrence of
the pattern $P$ if $p_i > p_j$ implies $a_{k_i} > a_{k_j}$
and we say that $a_{k_i}$ acts as $p_i$ in the occurrence, for all $i,j\in [m]$.
A generalized pattern is written as a classical pattern with a dash written between some of the
letters. If the letters $p_i$ and $p_{i+i}$ are not separated by a dash then we must also have
$k_{i+1} = k_i+1$ in an occurrence, or in other words,
the elements in the permutation acting as $p_i$ and $p_{i+i}$ must be adjacent.
We say that a pattern is contiguous if it has no dashes. We also say that a segment in a pattern is
contiguous if it has no dashes.
The permutation $\pi$ is said to avoid the pattern $P$ if $\pi$ has no occurrence of $P$.

\textbf{Example:} The subsequence 314 in the permutation 3124 is an occurrence of the pattern 21-3 since
314 is in same relative order as 213 and 3 and 1 are adjacent in 3124. The segment 21 in the pattern 21-3 is contiguous.

Knuth \cite{knuth} showed that a permutation can be sorted in one pass through a stack if and only if it avoids
the classical pattern 231 (2-3-1 in generalized pattern notation). He also showed that for each of
the patterns of length three, the permutations avoiding that pattern are counted by the Catalan numbers.

A pattern can be regarded as a function from $\mathcal{S}_n$ to $\mathbb{N}$ which counts
the number of occurrences of that pattern in a given permutation.
The pattern 2-1 in this sense is simply the number of inversions in a permutation and it is a Mahonian statistic.

Babson and Steingr\'{\i}msson \cite{bab} introduced generalized patterns
and showed that most Mahonian statistics in the literature can be written as a linear combination of generalized pattern functions.

Kitaev \cite{kita} introduced partially ordered generalized patterns when the pattern is built from
a permutation of a multiset instead of a permutation of a set.
Although one can write all partially ordered generalized patterns as a linear combination of generalized patterns
this notation turns out to very useful. For instance we can write (3-121-3) instead of (4-132-5) + (4-231-5) + (5-132-4) + (5-231-4).

Kitaev \cite{kita} studied  generalized patterns of the form $k$-$\sigma$, where $\sigma$ is contiguous and $k$ is the
largest letter in the pattern. He showed that the EGF (exponential generating function) for permutations avoiding the
pattern $k$-$\sigma$ can be written in terms of the EGF for permutations avoiding the pattern $\sigma$.

We will take this a bit further and show that the EGF for permutations avoiding a pattern of the form
$k$-$\sigma$-$k$ can also be written
in terms of the EGF of permutations avoiding the pattern $\sigma$.

A partition of a set $S$ is a collection of disjoint subsets of $S$, whose union is $S$.
The disjoint subsets are called blocks.
A bicolored set partition is a set partitions where each block has either of two colors assigned to it.

We will take a closer look at the patterns 3-12-3 and 3-121-3. We map avoiders of 3-12-3
bijectively to bicolored set partitions and avoiders of 3-121-3 bijectively to Dowling lattices.

It is well known that avoiding the pattern $12$ is equivalent to avoiding the pattern $1$-$2$.
Claesson \cite{clea} showed that avoiding the pattern $2$-$1$-$3$ is equivalent to avoiding
the pattern $2$-$13$. For partially ordered generalized patterns we will give a complete classification of
patterns $P$ and $Q$ such that avoiding the pattern $P$ is equivalent to avoiding the pattern $Q$.

\section{Preliminaries}
\textbf{Permutations.}  \ We define $\mathcal{S}_n$ as the set of all permutations of length $n$ and
$\mathcal{S} = \bigcup_{n=0}^\infty \mathcal{S}_n$ as the set of all permutations.

For a permutation $\pi \in \mathcal{S}_n$ we let $\pi_{\ell}$ be the part of the permutation
that lies to the left of the largest element and let $\pi_r$ be the part of the permutation
that lies to the right of the largest element. Thus, we can write  $\pi = \pi_{\ell}n\pi_r.$

For the permutation $\pi=637529184$ we have $\pi_{\ell} = 63752$ and $\pi_r = 184$.

For a permutation $\pi = a_1a_2\dots a_n$, we say that $a_i$ is a left-to-right maximum of $\pi$ or
that $\pi$ has a left-to-right maximum at $i$ if $a_j<a_i$ for all $j<i$.

For a permutation $\pi = a_1a_2\dots a_n$, we say that $a_i$ is a right-to-left maximum of $\pi$ or
that $\pi$ has a right-to-left maximum at $i$ if $a_j<a_i$ for all $j>i$.

For a permutation $\pi = a_1a_2\dots a_n$, we define $\LMAX(\pi)= \{a_i \, | \, a_i$ is a left-to-right maximum of $\pi\}$
and $\RMAX(\pi)= \{a_i \, | \, a_i$ is a right-to-left maximum of $\pi\}$

From these definitions we get for a permutation $\pi = \pi_{\ell}n\pi_r \in \mathcal{S}_n$:
\begin{displaymath} 
\LMAX(\pi) =  \LMAX(\pi_{\ell}n\pi_r)  = \LMAX(\pi_{\ell}) \cup \{n\} \\
\end{displaymath}
\begin{displaymath}
\RMAX(\pi) =  \RMAX(\pi_{\ell}n\pi_r)  = \RMAX(\pi_r) \cup \{n\}.
\end{displaymath}

\textbf{Example:} For the permutation $\pi=637529184$ we have:
\begin{displaymath}
\LMAX(637529184) = \LMAX(63752) \cup \{9\} = \LMAX(63) \cup \{7,9\} = \{6,7,9\} \\
\end{displaymath}
\begin{displaymath}
\RMAX(637529184) = \RMAX(184) \cup \{9\}  = \RMAX(4) \cup \{8,9\} = \{4,8,9\}.
\end{displaymath}

We let $[n] = \{1,2,\dots,n\}$ and define the function:
\begin{displaymath}
\Ins: \mathcal{S}_n \times [n+1] \times [n+1] \to \mathcal{S}_{n+1}
\end{displaymath}
\begin{displaymath}
\Ins(p_1p_2\dots p_n, i, j) \mapsto \widehat{p}_1 \widehat{p}_2 \dots \widehat{p}_{j-1} i \widehat{p}_j \dots \widehat{p}_n
\end{displaymath}
where $\widehat{p}_k = p_k$ if $p_k<i$ and $\widehat{p}_k = p_k+1$ otherwise. So we insert $i$ before $j$-th element in $\pi$ and increase all element greater or equal to $i$ by 1.

\textbf{Patterns.} \ A partially ordered generalized pattern $P$ is a permutation of
a multiset of ordered elements with a dash written between some of the letters in the permutation.
We call the permutation of the multiset the underlying permutation. We define the length of a pattern as 
the length of the underlying permutation.
Let $p_1p_2\dots p_m$ be the
underlying permutation of the pattern $P$ and $\pi = a_1a_2\dots a_n$ be a permutation.
For $1\le k_1<k_2<\dots<k_m \le n$ we say that $a_{k_1}a_{k_2}\dots a_{k_m}$ is an occurrence of
the pattern $P$ if $p_i > p_j$ implies $a_{k_i} > a_{k_j}$. Notice that if $p_i = p_j$ for some $i,j$ then
we have no restriction on the relative order of  $a_{k_i}$ and $a_{k_j}$.
If the letters $p_i$ and $p_{i+i}$ are not separated by a dash then we must also have
$k_{i+1} = k_i+1$ in an occurrence, or in other words,
the elements in the permutation acting as $p_i$ and $p_{i+i}$ must be adjacent.

For a permutation $\pi$ and pattern $P$ we say that $\pi$ forms $P$ if $\pi$ and $P$ have the same length
and $\pi$ has an occurrence of $P$. So an occurrence of a pattern $P$ is a subsequence of permutation that forms the pattern $P$.

For a pattern $P$ we define: 
\begin{displaymath}
\Av(P) = \{\pi\in\mathcal{S}  \mid \pi \; \mathrm{avoids} \; P\}
\end{displaymath}

\textbf{Exponential generating functions.} \ For a set $S$ and a weight function $\omega: \ S \to \mathbb{N}$
we define the exponential generating function (or $\EGF$) to be:
\begin{displaymath}
\sum_{s \in S}  \frac{x^{\omega(s)}}{{\omega(s)}!}
\end{displaymath}
The most natural weight function for permutations is the length and when we are considering the $\EGF$
for avoiders of a pattern $P$ we use the set $S=\Av(P)$.

For a pattern $P$ we denote by $\EGF_P$ the $\EGF$ for avoiders of the pattern $P$.

Now if a labeled structured having the $\EGF$ $h(x)$ is composed of two labeled structures having the $\EGF$'s $f(x)$ and $g(x)$
and the labeling of the two smaller labeled structures is independent then we have:
\begin{displaymath}
h(x) = f(x)g(x).
\end{displaymath}

Let $\mathbb{C}[[x]]$ be the ring of formal power series with coefficients in $\mathbb{C}$. 
We define the formal integral of $f(x)=\sum_{n=0}^\infty f_n x^n$ as follows:
\begin{displaymath}
\int f(x)dx = \sum_{n=0}^\infty f_n \frac{x^{n+1}}{n+1}.
\end{displaymath}

\section{The pattern $3$-$12$-$3$}

We will construct a bijection between permutations avoiding the
partially ordered generalized pattern $3$-$12$-$3$ and bicolored set partitions
(the parts have either of two colors).

We show that a permutation $\pi=\pi_{\ell} n \pi_r \in \mathcal{S}_n$ avoids $3$-$12$-$3$
if and only if $\pi_{\ell}$ avoids $3$-$12$ and  $\pi_r$ avoids $12$-$3$. We then use the fact
that ($3$-$12$)-avoiding permutations and ($12$-$3$)-avoiding permutations both map bijectively to
set partitions. We can then construct such a bijection in the following way.
For a permutation $\pi=\pi_{\ell} n \pi_r \in \mathcal{S}_n$, we map the ($3$-$12$)-avoiding
permutation $\pi_{\ell}$ to the corresponding set partition. We also map the ($12$-$3$)-avoiding
permutation $\pi_{r}$ to the corresponding set partition. We then assign one color to the set partition
derived from $\pi_{\ell}$ and the other color to the set partition derived from $12$-$3$. Together these two 
set partitions then form a bicolored set partition, and the described mapping is bijective.

\subsection{Two Lemmas}

Although we are dealing with the patterns $3$-$12$, $12$-$3$ and $3$-$12$-$3$, most of our
work will translate directly to similar pattern families. By similar pattern families we are referring
to patterns of the form $k$-$\sigma$, $\sigma$-$k$ and $k$-$\sigma$-$k$, where $\sigma$ is a contiguous
pattern and $k$ is greater than all the letters in $\sigma$. We start by proving important identities  
for avoidance of such patterns.

The following two lemmas can be proved as easily for $k$-$\sigma$-$k$ as for $3$-$12$-$3$ so we will prove
them for the more general version of the pattern. 

We require $\sigma$ to be contiguous and $k$ to be greater than all the letters of $\sigma$.

\begin{Lemma}\label{Lemma:gen1}

A permutation $\pi=\pi_{\ell} n \pi_r \in \mathcal{S}_n$ avoids $k$-$\sigma$ if and only if
$\pi_{\ell}$ avoids $k$-$\sigma$ and  $\pi_r$ avoids $\sigma$.

\textit{Proof:} (i) If $\pi_{\ell}$ has an occurrence of $k$-$\sigma$ then $\pi$ has an occurrence of
$k$-$\sigma$.
If $\pi_r$ has an occurrence of $\sigma$ then $n$ together with the elements forming the pattern $\sigma$
form an occurrence of $k$-$\sigma$ in $\pi$.

(ii) If $\pi$ has an occurrence of $k$-$\sigma$, then since all the letters in $\sigma$ are less than $k$,
$n$ can not act as any letter in $\sigma$. Also, 
since $\sigma$ is contiguous either all the elements acting as letters in $\sigma$ lie to the left of $n$ or
they all lie to the right of $n$.
In the first case the elements forming $k$-$\sigma$ in $\pi$ are all in $\pi_{\ell}$ and therefore form
$k$-$\sigma$ in $\pi_{\ell}$, in the latter case 
the elements forming the $\sigma$ part of $k$-$\sigma$
in $\pi$ form $\sigma$ in $\pi_r$. $\qed$

\end{Lemma}

\begin{Lemma}\label{Lemma:gen2}

A permutation $\pi=\pi_{\ell} n \pi_r \in \mathcal{S}_n$ avoids $k$-$\sigma$-$k$ if and only if
$\pi_{\ell}$ avoids $k$-$\sigma$ and  $\pi_r$ avoids $\sigma$-$k$.

\textit{Proof:} (i) If $\pi_{\ell}$ has an occurrence of $k$-$\sigma$ then the elements forming $k$-$\sigma$
together with $n$ form the pattern $k$-$\sigma$-$k$ in $\pi$ and $\pi$ has an occurrence of $k$-$\sigma$-$k$.
If $\pi_{r}$ has an occurrence of $\sigma$-$k$ then the elements forming $\sigma$-$k$
together with $n$ form the pattern $k$-$\sigma$-$k$ in $\pi$ and $\pi$ has an occurrence of $k$-$\sigma$-$k$.

(ii) If $\pi$ has an occurrence of $k$-$\sigma$-$k$, then since all the letters in $\sigma$ are less than $k$,
$n$ can not act as any letter in $\sigma$. Also, 
since $\sigma$ is contiguous either all the elements acting as letters in $\sigma$ lie to the left of $n$ or
they all lie to the right of $n$.
In the first case the elements forming the $k$-$\sigma$ part of $k$-$\sigma$-$k$ in $\pi$ are all 
in $\pi_{\ell}$ and therefore form $k$-$\sigma$ in $\pi_{\ell}$ and in the latter case
the elements forming the $\sigma$-$k$ part of $k$-$\sigma$-$k$ in $\pi$ are all 
in $\pi_r$ and therefore form $\sigma$-$k$ in $\pi_r$. $\qed$

\end{Lemma}

By using Lemma \ref{Lemma:gen1} recursively we get:

\begin{coro}\label{coro:gen1}

Let $\pi = M_1\pi_1M_2\pi_2\dots M_i\pi_i$, where $\{M_1, M_2, \dots, M_i\} = \LMAX(\pi)$ and 
$\pi_1,\pi_2,\dots,\pi_i$ are the blocks between left-to-right maxima. Then 
$\pi$ avoids $k$-$\sigma$ if and only if $\pi_j$ avoids $\sigma$, for $1 \le j \le i$. $\qed$

\end{coro}

\subsection{The patterns $3$-$12$ and $12$-$3$}

Let $\pi = M_1\pi_1M_2\pi_2\dots M_i\pi_i$, where $\{M_1, M_2, \dots, M_i\} = \LMAX(\pi)$ and
$\pi_1,\pi_2,\dots,\pi_i$ are the blocks between left-to-right maxima. Then Corollary \ref{coro:gen1}
gives that $\pi$ avoids $3$-$12$ if and only if $\pi_j$ avoids $12$, for $1 \le j \le i$.
That is equivalent to $\pi_j$ being decreasing, for $1 \le j \le i$.
Since $M_j \in \LMAX(\pi)$ we have that $M_j$ is larger than the letters of $\pi_j$ for $1 \le j \le i$
and $M_1 < M_2 < \dots < M_i$.

\begin{coro}\label{coro:set}

The permutation $\pi$ avoids $3$-$12$ if and only if $\pi$ is a sequence of decreasing sequences where the
first element in each sequence is greater than the first element in the preceding sequence. $\qed$
\end{coro}

Claesson \cite{clea} described a bijection between set partitions and ($3$-$12$)-avoiding permutations.

A set partition is mapped to a ($3$-$12$)-avoiding permutation as follows.
The elements within each block are written in a decreasing order and the blocks are written in
increasing order of the greatest element. 
From Corollary \ref{coro:set} we see that the resulting permutation avoids $3$-$12$. 

This process is easily reversible. For a given ($3$-$12$)-avoiding permutation
we use the left-to-right maxima and all the following elements until meeting the
next left-to-right maxima to construct the blocks in the set partition.

\textbf{Example:} A set partition and the corresponding ($3$-$12$)-avoiding permutation.
\begin{displaymath}
\{6, 3, 1\}\;  \{8,7\}\; \{9, 5 ,4 ,2\} \\
\end{displaymath}
\begin{displaymath}
6\,\; 3\,\; 1\,\; 8\,\; 7\,\; 9\,\; 5\,\; 4\,\; 2
\end{displaymath}

The number of blocks in the set partition will be equal to the number of left-to-right maxima
in the corresponding ($3$-$12$)-avoiding permutation.
This is because the largest element in each block will become a left-to-right maximum
and no other element can become a left-to-right maximum.

Claesson \cite{clea} also described a bijection between set partitions and ($12$-$3$)-avoiding permutations.
Instead of constructing this bijection, we construct a bijection between ($k$-$\sigma$)-avoiding permutations
and ($\sigma$-$k$)-avoiding permutations. We can then map ($12$-$3$)-avoiding
permutations bijectively to ($3$-$12$)-avoiding permutations which in turn can be
mapped bijectively to set partitions. As a consequence of this, we obtain Claesson's original bijection.

\begin{defi}
Let $\psi:$ $\mathcal{S} \to \mathcal{S}$ be the mapping:
\begin{eqnarray*}
\psi(\emptyset) & = & \emptyset  \\
\psi(\pi_{\ell} n \pi_r) & = & \pi_r \, n \, \psi(\pi_{\ell}).
\end{eqnarray*}
\end{defi}

\begin{Lemma}\label{Lemma:psi}
The map $\psi$ has the following properties:

\begin{itemize}
\item[(i)] $\psi$ preserves length.

\item[(ii)] $\psi$ is a bijection.

\item[(iii)] Let $\pi = M_1\pi_1M_2\pi_2\dots M_i\pi_i$, where $\{M_1, M_2, \dots, M_i\} = \LMAX(\pi)$ and
$\pi_1,\pi_2,\dots,\pi_i$ are the blocks between left-to-right maxima. Then,
\begin{displaymath}
\psi(M_1\pi_1M_2\pi_2\dots M_i\pi_i) =  \pi_iM_i\dots \pi_2M_2\pi_1M_1
\end{displaymath}

\item[(iv)] $\psi$ maps left-to-right maxima in $\pi$ to right-to-left maxima in $\phi(\pi)$, that is $\LMAX(\pi) = \RMAX(\psi(\pi))$.

\item[(v)] A permutation $\pi$ avoids the pattern $k$-$\sigma$ if and only if $\psi(\pi)$ avoids $\sigma$-$k$,
where $\sigma$ is contiguous pattern and $k$ is greater than the letters in $\sigma$.
\end{itemize}

\textit{Proof:}

(i) Now, $\psi(\pi)$ is a reordering of $\pi$ for all $\pi\in\mathcal{S}$ and therefor $\psi$ preserves length.

(ii) We show by induction on $n$ that $\psi$ is injective when it is restricted to
$\mathcal{S}_n$.

By definition $\psi$ is injective for $n=0$.

We assume that $\psi$ is injective for all $n<N$ for some $N\in\mathbb{N}$.

Let $\pi , \tau \in \mathcal{S}_N$ be such that $\psi(\pi) = \psi(\tau)$.
We write $\pi=\pi_{\ell}N\pi_r$ and $\tau=\tau_{\ell}N\tau_r$. Now
\begin{displaymath}
\pi_rN\psi(\pi_{\ell})= \psi(\pi_{\ell}N\pi_r)= \psi(\pi) = \psi(\tau)=\psi(\tau_{\ell}N\tau_r) =\tau_rN\psi(\tau_{\ell})
\end{displaymath}
so $\pi_r=\tau_r$ and $\psi(\pi_{\ell})=\psi(\tau_{\ell})$.
Since $|\pi_{\ell}|<N$ we have $\pi_{\ell}=\tau_{\ell}$.
Therefore $\pi = \tau$ and $\psi$ is injective for $n=N$.

We have shown that $\psi$ is injective on $\mathcal{S}_n$ for all $n\in\mathbb{N}$ and
(i) gives $\psi(\mathcal{S}_n)\subset \mathcal{S}_n$. Also $\mathcal{S}_n$ is finite so $\psi$
being injective on $\mathcal{S}_n$
is equivalent to $\psi$ being bijective on $\mathcal{S}_n$. We have therefore shown that $\psi$ is bijective.

(iii) We show by induction on the number of left-to-right maxima $i$ that
\begin{displaymath}
\psi(M_1\pi_1M_2\pi_2\dots M_i\pi_i) =  \pi_iM_i\dots \pi_2M_2\pi_1M_1.
\end{displaymath}
The statement obviously holds for $i=0$.

We assume that the statement holds for all $\tau \in \mathcal{S} $ such that $\#\LMAX(\tau)=j$ and
let $\pi = M_1\pi_1M_2\pi_2\dots M_j\pi_jM_{j+1}\pi_{j+1}$, then
$\pi_{\ell} = M_1\pi_1M_2\pi_2\dots M_j\pi_j$ so $\pi_{\ell}$ has $j$ left-to-right maxima and thus

\begin{displaymath}
\psi(\pi_{\ell}) = \psi(M_1\pi_1M_2\pi_2\dots M_j\pi_j) = \pi_jM_j\dots \pi_2M_2\pi_1M_1.
\end{displaymath}

We therefore have

\begin{displaymath}
\psi(\pi) = \psi(\pi_{\ell}M_{j+1}\pi_{j+1}) = \pi_{j+1}M_{j+1}\psi(\pi_{\ell}) = \pi_{j+1}M_{j+1}\pi_jM_j\dots \pi_2M_2\pi_1M_1.
\end{displaymath}

(iv) We let $\pi = M_1\pi_1M_2\pi_2\dots M_i\pi_i$.
From (iii) we have that $\psi(\pi) = \psi(M_1\pi_1M_2\pi_2\dots M_i\pi_i) =  \pi_iM_i\dots \pi_2M_2\pi_1M_1$.

We have $\LMAX(\pi) = \{M_1,M_2,\dots,M_i\}$ so all the letters in $\pi_j$ are less than $M_j$ for $1 \le j \le i$.

Therefore we have that no letter in $\pi_j$ can be in $\RMAX(\psi(\pi))$ and thus
\begin{displaymath}
\RMAX(\psi(\pi)) = \{M_1,M_2,\dots,M_i\}.
\end{displaymath}

(v) We show by induction on $n$ that $\pi \in \Av(k$-$\sigma)$ if and only if $\psi(\pi) \in \Av(\sigma$-$k)$,
for all $\pi \in \mathcal{S}_n$.

The statement obviously holds for $n=0$ since $\emptyset$ avoids all nontrivial patterns.

We assume that statement holds for all $n<N$ for some $N\in\mathbb{N}$.
Let $\pi = \pi_{\ell}N\pi_r \in \mathcal{S}_N$. From Lemma \ref{Lemma:gen1} we get the following:
\begin{displaymath}
\pi \in \Av(k\textrm{-}\sigma)  \Leftrightarrow \pi_{\ell} \in \Av(k\textrm{-}\sigma) \;\; \mathrm{ and } \;\; \pi_{r} \in \Av(\sigma)
\end{displaymath}
\begin{displaymath}
\psi(\pi) \in \Av(\sigma\textrm{-}k)  \Leftrightarrow \pi_r \in \Av(\sigma) \;\; \mathrm{ and } \;\; \psi(\pi_{\ell}) \in \Av(\sigma\textrm{-}k).
\end{displaymath}

But $|\pi_{\ell}| < N $ so $\pi_{\ell} \in \Av(k$-$\sigma) \Leftrightarrow \psi(\pi_{\ell}) \in \Av(\sigma$-$k) $
and therefore  $\pi \in \Av(k$-$\sigma) \Leftrightarrow \psi(\pi) \in \Av(\sigma$-$k)$.

So $\psi$ maps ($k$-$\sigma$)-avoiding permutations to ($\sigma$-$k$)-avoiding permutations. $\qed$

\end{Lemma}

Now (ii) together with (v) give that $\psi$ maps ($3$-$12$)-avoiding permutations bijectively to 
($12$-$3$)-avoiding permutations. If we look at the bijection between set partitions and ($3$-$12$)-avoiding 
permutations, described earlier, and its composition with $\psi$, we get a bijective mapping 
between set partitions and ($12$-$3$)-avoiding permutations.
From (iv) we see that the number of blocks in a set partition will be equal to the number of right-to-left maxima
in the resulting ($12$-$3$)-avoiding permutation. From (iii) we can see that this bijection maps
a set partitions to a ($12$-$3$)-avoiding permutation in the following way.
The largest element in each block is preceded by all the remaining elements in the block in a decreasing order, 
from left to right. The blocks are written in decreasing order of greatest letter.

\textbf{Example:} A set partition and the corresponding ($12$-$3$)-avoiding permutation.
    \begin{displaymath}
        \{5 ,4 ,2, 9\} \;  \{7,8\}\;  \{3, 1, 6\} \\
    \end{displaymath}
    \begin{displaymath}
        5\,\,\, 4\,\,\, 2\,\,\, 9\,\,\, 7\,\,\, 8\,\,\, 3\,\,\, 1\,\,\, 6
    \end{displaymath}

This bijection is due to Claesson \cite{clea}.

\subsection{The bijection}

From Lemma \ref{Lemma:gen2} we have that a permutation $\pi = \pi_{\ell}N\pi_r$  avoids $3$-$12$-$3$
if and only if $\pi_{\ell}$ avoids $3$-$12$ and  $\pi_r$ avoids $12$-$3$.
So the largest element in an ($3$-$12$-$3$)-avoiding permutation $\pi$ splits
the permutation into a ($3$-$12$)-avoiding permutation $\pi_{\ell}$ and a ($12$-$3$)-avoiding permutation $\pi_r$.

We can now map $\pi_{\ell}$ and $\pi_r$ each to a set partition using the bijections described earlier.
We then consider the two set partitions as a one single set partition but to distinguish between blocks
coming from $\pi_{\ell}$ and $\pi_r$, we assign one color to the blocks coming from $\pi_{\ell}$
and another color to the blocks coming from $\pi_r$. As a result we have a bicolored set partition and this
mapping is a bijection.

Furthermore the number of left-to-right maxima in $\pi_{\ell}$ is equal to
the number of parts of the first color and
the number of right-to-left maxima in $\pi_r$ is equal to
the number of parts of the second color.
But left-to-right maxima in $\pi$ are precisely the left-to-right maxima
in $\pi_{\ell}$ together with $n$ and right-to-left maxima in $\pi$ precisely right-to-left maxima in $\pi_r$
together with $n$.
So for a bicolored set partition and the corresponding ($3$-$12$-$3$)-avoiding permutation we have:
    \begin{displaymath}
        \textrm{The number of blocks of the first color } = \#\LMAX -1
    \end{displaymath}
    \begin{displaymath}
        \textrm{The number of blocks of the second color } = \#\RMAX -1.
    \end{displaymath}

\textbf{Example:} A ($3$-$12$-$3$)-avoiding permutation and the corresponding bicolored set partition.
    \begin{displaymath}
     \underbrace{ 6 \,\,\, 1 \,\,\,  9 \,\,\,  4 \,\,\, 2}_{\pi_{\ell}} \,\,\, 12 \,\,\, \underbrace{3 \,\,\, 11 \,\,\,8 \,\,\, 5 \,\,\, 10 \,\,\, 7}_{\pi_r}
    \end{displaymath}
    \begin{displaymath}
      \{6,1\}\;  \{9, 4, 2\} \quad \mathbf{\{3, 11\} \; \{8,5,10\} \; \{7\}}
    \end{displaymath}
(The blocks derived from $\pi_r$ in bold signify that they do not have the same color as the other blocks.)

\begin{coro}

Permutations in $\mathcal{S}_{n+1}$ avoiding $3$-$12$-$3$ are in one-to-one correspondence with bicolored
set partitions on the set $[n]$.

\end{coro}

\subsection{The generating functions}

\begin{defi}
Let $f(x)$, $g(x)$ and $h(x)$ be the exponential generating functions of 
$12$, ($3$-$12$) and ($3$-$12$-$3$)-avoiding permutations respectively. 
We denote by $f_n$, $g_n$ and $h_n$ the coefficients of $\frac{x^n}{n!}$ in $f(x)$, $g(x)$ and $h(x)$ respectively.
\end{defi}

There is only one permutation of length $n$ that avoids $12$, namely the decreasing permutation
so we have $f_n=1$ for all $n\in \mathbb{N}$ and thus $f(x) = e^x$.

We also note that permutations avoiding $3$-$12$ map bijectively to set partitions
and set partitions are counted by the Bell numbers, so $g_n$ is the $n$-th Bell number.

Lemma \ref{Lemma:gen1} gives that a permutation $\pi=\pi_{\ell} n \pi_r \in \mathcal{S}_n$ avoids $3$-$12$ if and only if $\pi_{\ell}$
avoids $3$-$12$ and  $\pi_r$ avoids $12$.
So each nonempty ($3$-$12$)-avoiding permutation is composed of a ($3$-$12$)-avoiding permutation, its
largest element and an $12$-avoiding permutation.
The largest letter $n$ functions as a separator between the ($3$-$12$)-avoiding permutation and
the $12$-avoiding permutation. The ($3$-$12$)-avoiding permutation $\pi_\ell$ and
$12$-avoiding permutation $\pi_r$ are independent and their combined size is $n-1$.
So we have that $g(x)f(x)$ has the same coefficients as $g(x)$ only shifted by one.

\begin{coro}

For all $n\ge0$,

\begin{displaymath}
g_{n+1} = \sum_{i=0}^n \binom{n}{i}g_if_{n-i} = \sum_{i=0}^n \binom{n}{i}g_i.
\end{displaymath}

This is the well known recursion for the Bell numbers.

We also get the integral equation:

\begin{displaymath}
  g(x) = 1+\int g(x)f(x)dx
\end{displaymath}

which has the unique solution:

\begin{displaymath}
  g(x) = \exp \int f(x)dx = \exp(e^x -1).
\end{displaymath}

\end{coro}

From Lemma \ref{Lemma:gen2} we have that each nonempty ($3$-$12$-$3$)-avoiding permutation
is composed of a ($3$-$12$)-avoiding permutation, its largest element and an ($12$-$3$)-avoiding permutation.
The largest letter $n$ acts as a separator between the ($3$-$12$)-avoiding permutation and
the ($12$-$3$)-avoiding permutation. The ($3$-$12$)-avoiding permutation $\pi_\ell$ and
($12$-$3$)-avoiding permutation $\pi_r$ are independent and their combined size is $n-1$.
So we have that $g(x)g(x)$ has the same coefficients as $h(x)$ only shifted by one.

\begin{coro}

For all $n\ge 0$:

\begin{displaymath}
h_{n+1} = \sum_{i=0}^n \binom{n}{i}g_ig_{n-i}.
\end{displaymath}

We also get the integral equation:

\begin{displaymath}
  h(x) = 1+\int g(x)^2dx = 1+\int \exp(2(e^x -1))dx.
\end{displaymath}

\end{coro}

\subsection{Avoiding the patterns $3$-$12$ and $k$-$(k-1)\dots21$}

Let $\pi = M_1\pi_1M_2\pi_2\dots M_i\pi_i$, where $\{M_1, M_2, \dots, M_i\} = \LMAX(\pi)$ and
$\pi_1,\pi_2,\dots,\pi_i$ are the blocks between left-to-right maxima. Then by using $12$
and $(k-1)\dots21$ as $\sigma$ Corollary \ref{coro:gen1}
gives that $\pi$ avoids the patterns $3$-$12$ and $k$-$(k-1)\dots21$
if and only if $\pi_j$ avoids the patterns $12$ and $(k-1)\dots21$,
for $1 \le j \le i$. That is equivalent to $\pi_j$ being decreasing of length less than $k-1$,
for $1 \le j \le i$.
Since $M_j \in \LMAX(\pi)$ we have that $M_j$ is larger than the letters of $\pi_j$ for $1 \le j \le i$
and $M_1 < M_2 < \dots < M_i$.

\begin{coro}

A permutation $\pi$ avoids $3$-$12$ and $k$-$(k-1)\dots21$ if and only if $\pi$ is a 
sequence of decreasing subsequences, each of length less than $k$ where the
first element in each sequence is greater than the first element in the preceding decreasing subsequence. $\qed$
\end{coro}

In the bijection between set partitions and ($3$-$12$)-avoiding 
permutations we described earlier, each decreasing sequence
is built from a block in the set partition. So restricting the length
of the decreasing sequence in a ($3$-$12$)-avoiding permutation is equivalent to restricting the
block size in the set partition.

\begin{coro}
Permutations that avoid the patterns $3$-$12$ and $k$-$(k-1)\dots21$ are in one-to-one correspondence with
set partitions with block size less than $k$. $\qed$
\end{coro}

From Lemma \ref{Lemma:psi}, using $\sigma = 12$ and $\sigma = (k-1)\dots21$, we get that
permutations avoiding the patterns $12$-$3$ and $(k-1)\dots21$-$k$ are in one-to-one correspondence to
permutations avoiding the patterns $3$-$12$ and $k$-$(k-1)\dots21$.

From Lemma \ref{Lemma:gen2} we have that a permutation $\pi = \pi_{\ell}N\pi_r$  avoids $3$-$12$-$3$
and $k$-$(k-1)\dots21$-$k$
if and only if $\pi_{\ell}$ avoids the patterns $3$-$12$ and $k$-$(k-1)\dots21$ and
$\pi_r$ avoids the patterns $12$-$3$ and $(k-1)\dots21$-$k$.

We can now map $\pi_{\ell}$ and $\pi_r$ each to a set partition with block size less than $k$
using the bijections described earlier.
We then consider the two set partitions as a single set partition but to distinguish between blocks
coming from $\pi_{\ell}$ and $\pi_r$, we assign one color to the blocks coming from $\pi_{\ell}$
and another color to the blocks coming from $\pi_r$. As a result we have bicolored set partitions
with block size less than $k$ and this mapping is a bijection.

\begin{coro}

Permutations in $\mathcal{S}_{n+1}$ avoiding $3$-$12$-$3$ and $k$-$(k-1)\dots21$-$k$ are in one-to-one correspondence with bicolored
set partitions with block size less than $k$. $\qed$

\end{coro}

\section{The pattern $k$-$\sigma$-$k$}

It turns out that we can use the same kind of arguments for patterns of the form $k$-$\sigma$-$k$
to obtain their exponential generating function from the exponential generating function of
permutations avoiding $\sigma$.

\subsection{The pattern $k$-$\sigma$}

We can use Lemma \ref{Lemma:gen1} on multiple $\sigma$'s. Also, avoiding a partially ordered pattern
is equivalent to avoiding multiple patterns.
Our results also applies to patterns of the form
$\sigma$-$k$, $1$-$\sigma$ and $\sigma$-$1$, because of symmetry.

\begin{coro} \

\begin{itemize}
\item[(i)] A permutation $\pi=\pi_{\ell} n \pi_r \in \mathcal{S}_n$ avoids the partially ordered patterns
$\{k_i$-$\sigma_i\}_{i\in I}$ if and only if $\pi_{\ell}$ avoids the partially ordered patterns
$\{k_i$-$\sigma_i\}_{i\in I}$ and  $\pi_r$ avoids the partially ordered patterns $\{\sigma_i\}_{i\in I}$.

\item[(ii)] A permutation $\pi=\pi_{\ell} n \pi_r \in \mathcal{S}_n$ avoids the partially ordered patterns
$\{\sigma_i$-$k_i\}_{i\in I}$ if and only if $\pi_{\ell}$ avoids the partially ordered patterns
$\{\sigma_i$-$k_i\}_{i\in I}$ and  $\pi_r$ avoids the partially ordered patterns $\{\sigma_i\}_{i\in I}$.

\item[(iii)] A permutation $\pi=\widehat{\pi}_{\ell} 1 \widehat{\pi}_r \in \mathcal{S}_n$ avoids the partially ordered patterns
$\{1$-$\sigma_i\}_{i\in I}$ if and only if $\widehat{\pi}_{\ell}$ avoids the partially ordered patterns
$\{1$-$\sigma_i\}_{i\in I}$ and  $\widehat{\pi}_r$ avoids the partially ordered patterns $\{\sigma_i\}_{i\in I}$.

\item[(iv)] A permutation $\pi=\widehat{\pi}_{\ell} 1 \widehat{\pi}_r \in \mathcal{S}_n$ avoids the partially ordered patterns
$\{\sigma_i$-$1\}_{i\in I}$ if and only if $\widehat{\pi}_{\ell}$ avoids the partially ordered patterns
$\{\sigma_i$-$1\}_{i\in I}$ and  $\widehat{\pi}_r$ avoids the partially ordered patterns $\{\sigma_i\}_{i\in I}$.
\end{itemize}

Where $\widehat{\pi}_{\ell}$ is the part of $\pi$ that lies to the left of 1 and
$\widehat{\pi}_r$ is the part of $\pi$ that lies to the right of 1.
\end{coro} \
Again by using Lemma \ref{Lemma:gen2} for multiple patterns and because of symmetry, we get the following
corollary.

\begin{coro} \

\begin{itemize}
\item[(i)] A permutation $\pi=\pi_{\ell} n \pi_r \in \mathcal{S}_n$ avoids the partially ordered patterns
$\{k_i$-$\sigma_i$-$k_i\}_{i\in I}$ if and only if $\pi_{\ell}$ avoids the partially ordered patterns
$\{k_i$-$\sigma_i\}_{i\in I}$ and  $\pi_r$ avoids the partially ordered patterns $\{\sigma_i$-$k_i\}_{i\in I}$.

\item[(ii)] A permutation $\pi=\widehat{\pi}_{\ell} 1 \widehat{\pi}_r \in \mathcal{S}_n$ avoids the partially ordered patterns
$\{1$-$\sigma_i$-$1\}_{i\in I}$ if and only if $\widehat{\pi}_{\ell}$ avoids the partially ordered patterns
$\{1$-$\sigma_i\}_{i\in I}$ and  $\widehat{\pi}_r$ avoids the partially ordered patterns $\{\sigma_i$-$1\}_{i\in I}$.

\end{itemize}

\end{coro}

\subsection{The generating functions}

\begin{defi}
Let $f(x)$, $g(x)$ and $h(x)$ be the exponential generating functions of $\sigma$, ($k$-$\sigma$) and 
($k$-$\sigma$-$k$)-avoiding permutations respectively. 
We denote by $f_n$, $g_n$ and $h_n$ the coefficients of $\frac{x^n}{n!}$ in $f(x)$, $g(x)$ and $h(x)$ respectively.
\end{defi}

We have from Lemma \ref{Lemma:gen1} that each nonempty ($k$-$\sigma$)-avoiding permutation is composed
of a ($k$-$\sigma$)-avoiding permutation, its largest element and a $\sigma$-avoiding permutation.
The largest letter $n$ functions as a separator between the ($k$-$\sigma$)-avoiding permutation and
the $\sigma$-avoiding permutation. The ($k$-$\sigma$)-avoiding permutation and
$\sigma$-avoiding permutation are independent and their combined size is $n-1$.
So we have that $g(x)f(x)$ has the same coefficients as $g(x)$ only shifted by one.

\begin{coro}\label{coro:egf1}

For all $n \ge 0$:

\begin{displaymath}
g_{n+1} = \sum_{i=0}^n \binom{n}{i}g_if_{n-i}.
\end{displaymath}

We also get the integral equation:

\begin{displaymath}
  g(x) = 1+\int g(x)f(x)dx
\end{displaymath}

which has the unique solution:

\begin{displaymath}
  g(x) = \exp\left(\int f(x)dx\right).
\end{displaymath}

\end{coro}

Again from Lemma \ref{Lemma:gen2} we have that each nonempty ($k$-$\sigma$-$k$)-avoiding permutation is composed
of a ($k$-$\sigma$)-avoiding permutation, its largest element and a ($\sigma$-$k$)-avoiding permutation.
The largest letter $n$ acts as a separator between the ($k$-$\sigma$)-avoiding permutation and
the ($\sigma$-$k$)-avoiding permutation. The ($k$-$\sigma$)-avoiding permutation and
($\sigma$-$k$)-avoiding permutation are independent and their combined size is $n-1$.
So we have that $g(x)^2$ has the same coefficients as $h(x)$, only shifted by one.

\begin{coro}\label{coro:egf2}

For all $n \ge 0$:

\begin{displaymath}
h_{n+1} = \sum_{i=0}^n \binom{n}{i}g_ig_{n-i}
\end{displaymath}

We also get the integral equation:

\begin{displaymath}
  h(x) = 1+\int g(x)^2dx.
\end{displaymath}

\end{coro}

\section{The pattern $3$-$121$-$3$}

For a group $G$ a $G$-labeled set $(S,\alpha)$ is a set $S$ together with a mapping $\alpha:S \to G$.
Two $G$-labeled sets $(S,\alpha)$ and $(S,\beta)$ are equivalent if there is $g\in G$ such that $\alpha = g \beta$.
We denote the equivalence class containing
$(S,\alpha)$ as $[S,\alpha]$. A partial $G$-partition of a set $S$ is
a set $\{[A_1,\alpha_1],[A_2,\alpha_2],\dots,[A_k,\alpha_k]\}$ where $A_i$ are disjoint subsets of $S$, for $i\in [k]$.

Dowling \cite{dowl} showed that partial $G$-partitions of a set $S$ form a lattice, for every group $G$ and set $S$. 
However, we do not need the partial ordering of the lattice just the number of elements in it. 
Also, the only group we are interested in is the two-element group $G=(\mathbb{Z}_2,+)$.
For a $G$-labeled set we often write
$\{a_{1_{\alpha(a_1)}},a_{2_{\alpha(a_2)}},\dots,a_{n_{\alpha(a_n)}}\}$ instead of
$(\{a_1,a_2,\dots,a_n\},\alpha)$.

\textbf{Example:} A partial $G$-partition of the set $[9]$.
\begin{displaymath}
[1_0, 6_1]\;  [7_0]\; [2_0, 4_0 ,5_1 ,9_0] 
\end{displaymath}
\begin{displaymath}
A_1=\{1,6\},\,\alpha_1(1)=0,\alpha_1(6)=1 \\
\end{displaymath}
\begin{displaymath}
A_2=\{7\},\,\alpha_2(7)=0 \\
\end{displaymath}
\begin{displaymath}
A_3=\{2,4,5,9\},\,\alpha_3(2)=0,\,\alpha_3(4)=0,\,\alpha_3(5)=1,\,\alpha_3(9)=0
\end{displaymath}

A $121$-avoiding permutation is a decreasing sequence followed by a increasing sequence and we
map them bijectively to the equivalence classes of $G$-labeled sets in the following way.
As a representative of the equivalence class we choose the labeling where the smallest element of the set
has the label $0$.
We call this labeling the standard labeling.
We then write the elements having the label $0$ down in a decreasing order followed by
the elements having the label $1$ in increasing order.

\textbf{Example:} The standard labeling of an equivalence class of $G$-labeled set and the corresponding $121$-avoiding permutation.
    \begin{displaymath}
        \{1_0,2_1,3_0,4_0,5_1,6_0,7_0,8_1,9_1\}
    \end{displaymath}
    \begin{displaymath}
        764312589
    \end{displaymath}

Let $\pi = M_1\pi_1M_2\pi_2\dots M_i\pi_i$, where $\{M_1, M_2, \dots, M_i\} = \LMAX(\pi)$ and
$\pi_1,\pi_2,\dots,\pi_i$ are the blocks between left-to-right maxima. Then Corollary \ref{coro:gen1}
gives that $\pi$ avoids $3$-$121$ if and only if $\pi_j$ avoids $121$, for $1 \le j \le i$.
Then $M_j\pi_j$ avoids $121$ as well, but then $M_j$ has the label 0 in the
corresponding standard labeling, 
since either $M_j$ is the smallest element in $M_j\pi_j$ or to the left of the smallest element.

Now, let $\pi = \pi_1M_1\pi_2M_2\dots \pi_iM_i$, where $\{M_1, M_2, \dots, M_i\} = \RMAX(\pi)$ and
$\pi_1,\pi_2,\dots,\pi_i$ are the blocks between right-to-left maxima. Then corollary \ref{coro:gen1}
and Lemma \ref{Lemma:psi}
give that $\pi$ avoids $121$-$3$ if and only if $\pi_j$ avoids $121$, for $1 \le j \le i$.
Then $\pi_jM_j$ avoids $121$ as well, but we only map $\pi_jM_j$ to the corresponding standard labeling
if $\pi_j$ is nonempty. Then $M_j$ has the label 1 in the
corresponding standard labeling, since $M_j$ is to the right of the smallest element.
The remaining $M_j$ we map to the set $\Phi=\{M_j: \pi_j$ is empty$\}$.

From Lemma \ref{Lemma:gen2} we have that a permutation $\pi = \pi_{\ell}N\pi_r$  avoids $3$-$121$-$3$
if and only if $\pi_{\ell}$ avoids $3$-$121$ and  $\pi_r$ avoids $121$-$3$.
So the largest element in a ($3$-$121$-$3$)-avoiding permutation $\pi$ splits
the permutation into a ($3$-$121$)-avoiding permutation $\pi_{\ell}$ and a ($121$-$3$)-avoiding permutation $\pi_r$.
Thus, we can map ($3$-$121$-$3$)-avoiding permutations bijectively to partial $G$-partitions.
The blocks having the label 0 on their greatest element, in their standard labeling, come
from $\pi_{\ell}$ and the blocks having the label 1 on their greatest
element, in their standard labeling, come from $\pi_r$ as described above. Also, in the set $\Phi$ we have the
elements who are not in the partial  $G$-partition.

\textbf{Example:} A ($3$-$121$-$3$)-avoiding permutation and the corresponding partial $G$-partition. $\Phi=\{3,8\}$
\begin{displaymath}
 7 \;\;\; 9\;4\;2\;5\;\;\;10\;\;\; 8 \;\;\;1\;6 \;\;\;3
\end{displaymath}
\begin{displaymath}
[7_0]\;   [2_0, 4_0 ,5_1 ,9_0]\; [1_0, 6_1]
\end{displaymath}

\begin{coro}

Permutations in $\mathcal{S}_{n+1}$ avoiding $3$-$121$-$3$ are in one-to-one correspondence with
partial $G$-partitions on the set $[n]$, where $G$ is the two-element group.

\end{coro}

\begin{coro}

The number of permutations in $\mathcal{S}_{n+1}$ avoiding $3$-$121$-$3$ is equal to
the number of elements in the Dowling lattice over the set $[n]$ generated by the two-element group.

\end{coro}

Now since $\EGF_1(x) = 1$ we get from Corollary \ref{coro:egf1}
\begin{displaymath}
\EGF_\textrm{2-1}(x) = \exp\left(\int \EGF_1(x) dx\right) = e^x.
\end{displaymath}
In next section we will prove that avoiding the pattern 212 is equivalent to 
to avoiding the pattern 2-1-2, so from Corollary \ref{coro:egf2} we get
\begin{displaymath}
\EGF_{121}(x) = \EGF_{212}(x) = \EGF_{\textrm{2-1-2}}(x) = 1+\int \EGF_{\textrm{2-1}}(x)^2dx = \frac{e^{2x}+1}{2}.
\end{displaymath}
We now use Corollary \ref{coro:egf1} with $\sigma = 121$ and obtain
\begin{displaymath}
\EGF_\textrm{3-121}(x) = \exp\left(\int \EGF_{121}(x) dx\right) = \exp\left(\frac{e^{2x}+2x-1}{4}\right)
\end{displaymath}
and thus from Corollary \ref{coro:egf2} we get
\begin{displaymath}
\EGF_\textrm{3-121-3}(x) = 1 + \int \exp\left(\frac{e^{2x}+2x-1}{2}\right)dx.
\end{displaymath}

\section{Avoidance of a single partially ordered pattern}
We are particularly interested in contiguous patterns like $12$.
Avoiding the pattern $12$ is equivalent to avoiding the pattern $1$-$2$.
Claesson \cite{clea} showed that avoiding the pattern $2$-$1$-$3$  is equivalent to avoiding
the pattern $2$-$13$. So a natural question to ask is for which pairs of different partially ordered
patterns $P$ and $Q$ it is equivalent to avoid $P$ and to avoid $Q$?

For partially ordered generalized patterns $P$ and $Q$ such that
to avoid $P$ is equivalent to avoid $Q$ we observe that $P'=Q'$,
where $P'$ and $Q'$ are obtained from $P$ and $Q$ by removing all dashes.

All patterns of fixed length consisting only of the letter 1 have the same avoiders. So patterns consisting
only of the letter 1 are trivial with respect to addition or removal of a single dash.

We will prove the following two propositions, which give a complete classification of when two
partially ordered patterns have the same avoiders.

For a set $A$ we denote all string over $A$ by $A*$.

\begin{prop} \label{prop:av}
Let $P=\sigma_1$-$\sigma_2$-$\sigma_3$-$\sigma_4$ and $Q=\sigma_1$-$\sigma_2\sigma_3$-$\sigma_4$ be
two nontrivial partially ordered generalized patterns, where $\sigma_2$ and $\sigma_3$ are contiguous.
By symmetry we can assume that
the rightmost letter in $\sigma_2$ is less or equal to the leftmost letter in $\sigma_3$.
Then $\Av(P) = \Av(Q)$ if and only if $P$ has one of the following properties:
\begin{itemize}
\item[(1)] $\sigma_2 = 1\dots1$, $\sigma_3 = 2$, where $P' \in \{1,2\}^*$.
\item[(2)] $\sigma_2 = 1$, $\sigma_3 = 2\dots2$, where $P' \in \{1,2\}^*$.
\item[(3)] $\sigma_2 = 2\dots21$, $\sigma_3 = 2\dots2$, where $P' \in \{1,2\}^*$ and 1 appears once in $P$.
\item[(4)] $\sigma_2 = 1\dots1$, $\sigma_3 = 21\dots1$, where $P' \in \{1,2\}^*$ and 2 appears once in $P$.
\item[(5)] $\sigma_2 = 1\dots1$, $\sigma_3 = 3$, where $P' \in \{1,2,3\}^*$ and 2 appears once in $P$.
\item[(6)] $\sigma_2 = 1$, $\sigma_3 = 3\dots3$, where $P' \in \{1,2,3\}^*$ and 2 appears once in $P$.
\end{itemize}
\end{prop}

\begin{prop} \label{prop:av2}
Let $P$ and $Q$ be nontrivial partially ordered generalized patterns such that $\Av(P) = \Av(Q)$.
Then $Q$ can be obtained from $P$ with a succession of additions and removals of single dashes
so that at each step of the way
the set of avoiders of the pattern before and after the addition or the removal of the single dash
are the same.

\end{prop}

\begin{coro}
The only generalized patterns $P$ and $Q$ such that $\Av(P) = \Av(Q)$
are given by the following two examples.
\end{coro}

\textbf{Example:} From (1) or (2) in Proposition \ref{prop:av}, we get:
    \begin{displaymath}
        \Av(\textrm{1-2}) = \Av(12).
    \end{displaymath}

\textbf{Example:} From (5) or (6) in Proposition \ref{prop:av}, we get:
    \begin{displaymath}
        \Av(\textrm{2-1-3}) = \Av(\textrm{2-13}).
    \end{displaymath}

\begin{coro}
There are infinitely many
partially ordered generalized patterns $P$ and $Q$ such that $\Av(P) = \Av(Q)$.
\end{coro}

\textbf{Example:} By using first (1) or (2) and then (3) in Proposition \ref{prop:av}, we get:
    \begin{displaymath}
        \Av(\textrm{2-1-2}) = \Av(\textrm{21-2}) = \Av(212).
    \end{displaymath}

\textbf{Example:} By using (1) or (2) in Proposition \ref{prop:av} repeatedly, we get:
    \begin{displaymath}
        \Av(\textrm{12-12}) = \Av(\textrm{1-2-12})= \Av(\textrm{1-2-1-2})= \Av(\textrm{1-21-2}).
    \end{displaymath}

Our proof of Proposition \ref{prop:av} is a case-by-case proof and is rather lengthy.
Therefore we skip it until the end of this section. But that proof gives the following lemma for free.

\begin{Lemma}\label{Lemma:pop}
Let $P=\sigma_1$-$\sigma_2$ and $Q=\sigma_1\sigma_2$ be
two partially ordered generalized patterns.
The patterns $P$ and $Q$ have the same avoiders if and only if for all permutations $AXB$
such that $AB$ forms $Q$, $AXB$ has occurrence of $Q$,
where $|A| = |\sigma_1|$ and $|B| = |\sigma_2|$ ($A$ forms $\sigma_1$ and $B$ forms $\sigma_2$).
\end{Lemma}

We can then use Lemma \ref{Lemma:pop} prove the following lemma.

\begin{Lemma}\label{Lemma:pop2}
Let $P=\sigma_1\sigma_2$, $Q=\sigma_3$-$\sigma_4$ and $R=\sigma_1$-$\sigma_2$
be three partially ordered generalized patterns,
where $P$ and $Q$ have the same underlying permutation, $\Av(P)=\Av(Q)$, $|\sigma_1|=|\sigma_3|$ and $|\sigma_2|=|\sigma_4|$.
Then $\Av(P)=\Av(R)$.

\textit{Proof:} We take a permutation $AXB$ where $AB$ forms $\sigma_1\sigma_2$, $|A| = |\sigma_1$| and $|B| = |\sigma_2|$
then $AB$ forms $\sigma_3\sigma_4$, $|A| = |\sigma_3$| and $|B| = |\sigma_4|$ so $AXB$ has an occurrence of
$Q=\sigma_3$-$\sigma_4$. Since $\Av(P)=\Av(Q)$, $AXB$ has an occurrence of $P=\sigma_1\sigma_2$.
This holds for every such permutation $AXB$ and thus we get from Lemma \ref{Lemma:pop} that $\Av(P)=\Av(R)$. $\qed$
\end{Lemma}

Proposition \ref{prop:av2} is in fact a corollary to Lemma \ref{Lemma:pop2}.
Let $P$ and $Q$ be two partially ordered patterns, so that $\Av(P)=\Av(Q)$. Let $R$ be the pattern
having the same underlying permutation as $P$ and $Q$. Let $R$ have dash where either $P$ or $Q$ has dash.
We can then obtain $R$ from either $P$ or $Q$ by adding one dash at a time and from Lemma \ref{Lemma:pop2}
we get that set of avoiders at each step is invariant.

We will now prove Proposition \ref{prop:av} in the two following lemmas.

\begin{Lemma}
Let $P=\sigma_1$-$\sigma_2$-$\sigma_3$-$\sigma_4$ and $Q=\sigma_1$-$\sigma_2\sigma_3$-$\sigma_4$ be
two nontrivial partially ordered generalized patterns, where $\sigma_2$ and $\sigma_3$ are contiguous.
We write $\sigma_2=a_i\dots a_2a_1$ and $\sigma_3=b_1b_2\dots b_j$.
By symmetry we can assume that
the rightmost letter in $\sigma_2$ is less than or equal to the leftmost letter in $\sigma_3$
(from now on we let $a_1 \le b_1$).
We let $n = |P|$ and $k=|\sigma_1$-$\sigma_2|$.
In the following cases $\Av(P) \ne \Av(Q)$.

\begin{itemize}
\item[(i)] There exists a letter $c$ in $P$ such that $c < a_1$ or such that $c > b_1$.
\item[(ii)] There exist two letters $c_1$ and $c_2$ in $P$ such that $a_1 < c_1,c_2 < b_1$.
\item[(iii)] There exists a letter $c$ in $\sigma_2$ such that $c > a_1$ and there exists
a letter in $P$ other than $a_1$ less than $b_1$. Or there exists a letter $c$ in $\sigma_3$
such that $c < b_1$ and there exists a letter in $P$ other than $b_1$ greater than $a_1$.
\item[(iv)] If $a_1 = a_2$ and $b_1 = b_2$.
\end{itemize}
\textit{Proof:}

(i) There exists a letter $c$ in $P$ such that $c < a_1$ or such that $c > b_1$.

Let $c$ be a letter in $P$ such that $c > b_1$ and let $\pi=p_1p_2\dots p_n$ be a permutation such that
$\pi \notin \Av(P)$. We look at $\tau =\Ins(\pi,n+1,k+1)=p_1p_2\dots p_k(n+1)p_{k+1}\dots p_n$.
Now if $\tau$ has an occurrence of $Q$ then either $p_k$ acts as $a_1$ and $n+1$ acts as $b_1$ or
$n+1$ acts as $a_1$ and $p_{k+1}$ acts as $b_1$. Now $c>a_1,b_2$ so $n+1$ can neither act as
$a_1$ nor $b_1$. Therefore $\tau$ does not have an occurrence of $Q$.

We can construct a counterexample in a similar way for the case where there exists a letter $c$ in $P$
such that $c < a_1$.

\textbf{Example:} $\Av(\textrm{3-2-1}) \ne \Av(\textrm{32-1})$ since 321 forms  32-1
but 4\textbf{1}32 has only occurrence of 3-2-1

\textbf{Example:} $\Av(\textrm{3-2-1}) \ne \Av(\textrm{3-21})$ since 321 forms  3-21
but 32\textbf{4}1 has only occurrence of 3-2-1

(ii) There exist two letters $c_1$ and $c_2$ in $P$ such that $a_1 < c_1,c_2 < b_1$.

Let $c_1$ and $c_2$ be two letters in $P$ such that $a_1 < c_1 , c_2 < b_1$ and let  $\pi=p_1p_2\dots p_n$
be a permutation such that $\pi \notin \Av(P)$ and
let $\gamma_1$ and $\gamma_2$ be the elements in $\pi$ that act as $c_1$ and $c_2$. We can assume that
$\gamma_1<\gamma_2$.
We let $X=\{a$ letter in $P$ $|$ $a \le a_1\}$ and $Y=\{a$ letter in $P$ $|$ $a \ge b_1\}$.
We look at $\tau =\Ins(\pi, \gamma_2,k+1)=\widehat{p}_1\widehat{p}_2\dots \widehat{p}_k\gamma_2\widehat{p}_{k+1}\dots \widehat{p}_n$, where
$\widehat{p}_i = p_i+1$ if $p_i \ge \gamma_2$ otherwise $\widehat{p}_i = p_i$.

All the elements in $\tau$, derived from the elements that act as a letter from $X$ and $\gamma_1$ in $\pi$,
are all less than $\gamma_2$. So there are at least $|X|+1$ elements less than $\gamma_2$ in $\tau$.

Also, all the elements in $\tau$, derived from the elements that act as a letter from $Y$ and $\gamma_2$ in $\pi$,
are all greater than $\gamma_2$. So there are at least $|Y|+1$ elements greater than $\gamma_2$ in $\tau$.

Now if $\tau$ has an occurrence of $Q$ then either $p_k$ acts as $a_1$ and $\gamma_2$ acts as $b_1$,
or $\gamma_2$ acts as $a_1$ and $p_{k+1}$ acts as $b_1$.
Also, all but one element in $\tau$ are used in the occurrence
so at least $|X|$ of the $|X|+1$ elements, that are less than $\gamma_2$, will be used in the occurrence
and at least $|Y|$ of the $|Y|+1$ elements, that are greater than $\gamma_2$, will be used in the occurrence.
But there are at least $|X|$ elements in $\tau$ less than $\gamma_2$ so they have to act as the letters in $X$
and therefore can $\gamma_2$ not act as $a_1\in X$.
Also, there are at least $|Y|$ elements in $\tau$ greater than $\gamma_2$ so they have to act as the letters in $Y$
and therefore can $\gamma_2$ not act as $b_1\in X$.

We have thus shown that $\tau$ does not have an occurrence of $Q$.

\textbf{Example:} $\Av(\textrm{1-3-2-2}) \ne \Av(\textrm{13-2-2})$ since 1432 forms 13-2-2
but 1\textbf{3}542 has only occurrence of 1-3-2-2

(iii) There exists a letter $c$ in $\sigma_2$ such that $c > a_1$ and there exists
a letter in $P$ other than $a_1$ less than $b_1$. Or there exists a letter $c$ in $\sigma_3$
such that $c < b_1$ and there exists a letter in $P$ other than $b_1$ greater than $a_1$.

Let $c$ be a letter in $\sigma_2$ such that $c > a_1$ and let us assume that there exists
a letter in $P$ other than $a_1$ less than $b_1$.
We can assume that $a_1 = \dots = a_i$ and $a_{i+1} > a_1$. We can also assume there are no letters
less than $a_1$ in $P$ because otherwise (i) gives $\Av(P) \ne \Av(Q)$.

Let  $\pi=p_1p_2\dots p_{k-i}1p_{k-i+2}\dots p_kp_{k+1}\dots p_n$ be a permutation such that $\pi \notin \Av(P)$.
Such a permutation exist because 1 acts as $a_i$ and there is no letter less than $a_i=a_1$ in $P$.
We look at $\tau =\Ins(\pi, 2,k+1)=\widehat{p}_1\widehat{p}_2\dots \widehat{p}_{k-i}1\widehat{p}_{k-i+2}\dots \widehat{p}_k2\widehat{p}_{k+1}\dots \widehat{p}_n$, where
$\widehat{p}_j = p_j+1$ if $p_j \ge 2$ otherwise $\widehat{p}_j = p_j$. Now if $\tau$ has an occurrence of $Q$ then either
\begin{itemize}
\item[]1 acts as $a_{i+1}$,
\item[]$\widehat{p}_{k-i+2}$ acts as $a_i$,
\item[]$\quad\quad\quad\vdots$
\item[]$\widehat{p}_k$ acts as $a_2$,
\item[]2 acts as $a_1$
\end{itemize}
or
\begin{itemize}
\item[]1 acts as $a_{i}$,
\item[]$\widehat{p}_{k-i+2}$ acts as $a_{i-1}$,
\item[]$\quad\quad\quad\vdots$
\item[]$\widehat{p}_k$ acts as $a_1$,
\item[]2 acts as $b_1$.
\end{itemize}
Now 1 can not act as $a_{i+1}$ because $a_{i+1}>a_1$ and 2 can not act as $b_1$ because there are two different elements less than $b_1$ in $Q$.

We have thus shown that $\tau$ does not have an occurrence of $Q$.

\textbf{Example:} $\Av(\textrm{21-3}) \ne \Av(\textrm{213})$ since 231 forms  231
but 31\textbf{2}4 has only occurrence of 21-3

\textbf{Example:} $\Av(\textrm{31-3-2}) \ne \Av(\textrm{313-2})$ since 4132 forms  313-2
but 51\textbf{2}43 has only occurrence of 31-3-2

If $c$ there exists a letter $c$ in $\sigma_3$
such that $c < b_1$ and there exists a letter in $P$ other than $b_1$ greater than $a_1$. We can create
a counterexample in the similar way.

\textbf{Example:} $\Av(\textrm{1-32}) \ne \Av(\textrm{132})$ since 132 forms  132
but 1\textbf{3}42 has only occurrence of 1-32

\textbf{Example:} $\Av(\textrm{1-31-2}) \ne \Av(\textrm{131-2})$ since 1423 forms  131-2
but 1\textbf{4}523 has only occurrence of 1-31-2

(iv) If $a_1 = a_2$ and $b_1 = b_2$.

Let  $\pi=p_1p_2\dots p_kp_{k+1}\dots p_n$ be a permutation such that $\pi \notin \Av(P)$.
We look at $\tau =\Ins(\Ins(\pi, 1,k+1) ,n+2,k+1)=\widehat{p}_1\widehat{p}_2\dots \widehat{p}_k (n+2)1 \widehat{p}_{k+1}\dots \widehat{p}_n$, where
$\widehat{p}_i = p_i+1$.
Now if $\tau$ has an occurrence of $Q$ then one of the following holds:
\begin{itemize}
\item[] $n+2$ acts as $a_2$ and $1$  acts as $a_1$ or
\item[] $n+2$ acts as $a_1$ and $1$  acts as $b_1$ or
\item[] $n+2$ acts as $b_1$ and $1$  acts as $b_2$
\end{itemize}
But $n+2$ can not act as $a_1$ or $a_2$ nor can 1 act as $b_1$ or $b_2$,
since $n+2$ is the largest element in $\tau$ and 1 is the smallest element in $\tau$.

We have thus shown that $\tau$ does not have an occurrence of $Q$.

\textbf{Example:} $\Av(\textrm{11-22}) \ne \Av(\textrm{1122})$ since 1234 forms  1122
but 23\textbf{61}45 has only occurrence of 11-22

$\qed$
\end{Lemma}

Now if $\Av(P) = \Av(Q)$ then (i) together with (ii) give that only the letters 1,2 and 3 can appear in $P$.
Furthermore if the letter 3 appears in $P$ then the letter 2 appears exactly once. From (iii) we get that if
there exist a letter greater than the rightmost letter in $\sigma_2$ then all letters in $P$ except the
rightmost letter in $\sigma_2$ are equal to the leftmost letter in $\sigma_3$. Also from (iii) we get that if
there exist a letter less than the leftmost letter in $\sigma_3$ then all letters in $P$ except the
leftmost letter in $\sigma_2$ are equal to the rightmost letter in $\sigma_3$.
Otherwise $\sigma_2$ and $\sigma_3$ each consists only of a single letter
and then from (iv) we get that only one of $\sigma_2$ and $\sigma_3$ can have length greater than 1.
So $P$ has one of the following properties:

Otherwise only one letter appears in $\sigma_2$ and only one letter appears in $\sigma_3$ and then we get
from (iv) we get that
only one of $\sigma_2$ and $\sigma_3$ can have length greater than 1. So $P$ has one of the following properties
\begin{itemize}
\item[(1)] $\sigma_2 = 1\dots1$, $\sigma_3 = 2$, where $P' \in \{1,2\}^*$.
\item[(2)] $\sigma_2 = 1$, $\sigma_3 = 2\dots2$, where $P' \in \{1,2\}^*$.
\item[(3)] $\sigma_2 = 2\dots21$, $\sigma_3 = 2\dots2$, where $P' \in \{1,2\}^*$ and 1 appears once in $P$.
\item[(4)] $\sigma_2 = 1\dots1$, $\sigma_3 = 21\dots1$, where $P' \in \{1,2\}^*$ and 2 appears once in $P$.
\item[(5)] $\sigma_2 = 1\dots1$, $\sigma_3 = 3$, where $P' \in \{1,2,3\}^*$ and 2 appears once in $P$.
\item[(6)] $\sigma_2 = 1$, $\sigma_3 = 3\dots3$, where $P' \in \{1,2,3\}^*$ and 2 appears once in $P$.
\end{itemize}

We only need to show that for a pattern $P$ having one of the previous properties $\Av(P)=\Av(Q).$

\begin{Lemma}
Let $P=\sigma_1$-$\sigma_2$-$\sigma_3$-$\sigma_4$ and $Q=\sigma_1$-$\sigma_2\sigma_3$-$\sigma_4$ be
two nontrivial partially ordered generalized patterns, where $\sigma_2$ and $\sigma_3$ are contiguous.
We write $\sigma_2=a_i\dots a_2a_1$ and $\sigma_3=b_1b_2\dots b_j$.
By symmetry we can assume that
the rightmost letter in $\sigma_2$ is less than or equal to the leftmost letter in $\sigma_3$
(from now on we let $a_1 \le b_1$).
We let $n = |P|$ and $k=|\sigma_1$-$\sigma_2|$.
In the following cases $\Av(P) = \Av(Q)$.

(i) If all the letters in $P$ except $a_1$ are equal and all those letters are greater than $a_1$.
(ii) If $|\sigma_2| = 1$, $b_i=b_1$ for all $i$ and there exist at most one single letter $c$ in $P$ so that
$a_1<c<b_1$ and for all other letters $d$ in $P$ either $d=a_1$ or $d=b_1$.

\textit{Proof:}

(i) If all the letters in $P$ except $a_1$ are equal and all those letters are greater than $a_1$.

Let $\pi=p_1p_2\dots p_N \notin \Av(P)$ and
\begin{displaymath}
q_1q_2\dots q_{|\sigma_1|}q_{|\sigma_1|+1}\dots q_{|\sigma_1|+|\sigma_2|}q_{|\sigma_1|+|\sigma_2|+1}\dots q_{|\sigma_1|+|\sigma_2|+|\sigma_3|}\dots q_{|P|}
\end{displaymath}
be an occurrence of $P$ in $\pi$. Where $q_{|\sigma_1|+1}\dots q_{|\sigma_1|+|\sigma_2|}$ and $q_{|\sigma_1|+|\sigma_2|+1}\dots q_{|\sigma_1|+|\sigma_2|+|\sigma_3|}$ are adjacent in $\pi$
and $q_{|\sigma_1|+|\sigma_2|}$ is the smallest of elements in the occurrence.
We let $\widehat{q}_{|\sigma_1|+|\sigma_2|}$ be the smallest element of the elements
$q_{|\sigma_1|+|\sigma_2|}$ and the elements between $q_{|\sigma_1|+|\sigma_2|}$ and $q_{|\sigma_1|+|\sigma_2|+1}$ in $\pi$.
We choose as $\widehat{q}_{|\sigma_1|+1}$,\dots, $\widehat{q}_{|\sigma_1|+|\sigma_2|-1}$ the elements immediately  to the left of $\widehat{q}_{|\sigma_1|+|\sigma_2|}$ in $\pi$
and as $\widehat{q}_{|\sigma_1|+|\sigma_2|+1}$,\dots,$\widehat{q}_{|\sigma_1|+|\sigma_2|+|\sigma_3|}$ we choose the elements immediately  to the left of $\widehat{q}_{|\sigma_1|+|\sigma_2|}$ in $\pi$.
For the remaining $i$ we set $\widehat{q_i}=q_i$. Then
\begin{displaymath}
\widehat{q}_1\widehat{q}_2\dots \widehat{q}_{|\sigma_1|}\widehat{q}_{|\sigma_1|+1}\dots \widehat{q}_{|\sigma_1|+|\sigma_2|}\widehat{q}_{|\sigma_1|+|\sigma_2|+1}\dots \widehat{q}_{|\sigma_1|+|\sigma_2|+|\sigma_3|}\dots \widehat{q}_{|P|}
\end{displaymath}
is an occurrence of $Q$ in $\pi$.

We have thus shown that $\Av(P)=\Av(Q)$ for patterns $P$ of the type (3) and
we can prove in a similar way that the same holds for patterns $P$ of the type (4).

(ii) If $|\sigma_2| = 1$, $b_i=b_1$ for all $i$ and there exist at most one single letter $c$ in $P$ so that
$a_1<c<b_1$ and for all other letters $d$ in $P$ either $d=a_1$ or $d=b_1$.

Let $\pi=p_1p_2\dots p_N \notin \Av(P)$ and
\begin{displaymath}
q_1q_2\dots q_{|\sigma_1|}q_{|\sigma_1|+1}q_{|\sigma_1|+|\sigma_2|+1}\dots q_{|\sigma_1|+|\sigma_2|+|\sigma_3|}\dots q_{|P|}
\end{displaymath}
be an occurrence of $P$ in $\pi$. We set $c$ as the letter in $P$ so that $a_1<c<b_1$ if such a letter exists.
Otherwise we set as $c$ as the minimal element in $\pi$ acting as a letter $l=b_1$ in $P$.
We let $\widehat{q}_{|\sigma_1|+1}$ be the rightmost element less than $c$ of the elements
$q_{|\sigma_1|+1}$ and the elements between $q_{|\sigma_1|+1}$ and $q_{|\sigma_1|+|\sigma_2|+1}=q_{|\sigma_1|+2}$ in $\pi$.
As $\widehat{q}_{|\sigma_1|+|\sigma_2|+1}$,\dots,$\widehat{q}_{|\sigma_1|+|\sigma_2|+|\sigma_3|}$ we choose the elements immediately  to the left of $\widehat{q}_{|\sigma_1|+1}$ in $\pi$.
For the remaining $i$ we set $\widehat{q_i}=q_i$. Then
\begin{displaymath}
\widehat{q}_1\widehat{q}_2\dots \widehat{q}_{|\sigma_1|}\widehat{q}_{|\sigma_1|+1}\widehat{q}_{|\sigma_1|+|\sigma_2|+1}\dots \widehat{q}_{|\sigma_1|+|\sigma_2|+|\sigma_3|}\dots \widehat{q}_{|P|}
\end{displaymath}
is an occurrence of $Q$ in $\pi$. $\qed$
\end{Lemma}

We have thus shown that $\Av(P)=\Av(Q)$ for patterns $P$ of the type (2) or (6) and
we can prove in a similar fashion that the same holds for patterns $P$ of the type (1) or (5).

In each case where $\Av(P) \ne \Av(Q)$ we have shown that there exists a permutation $AXB$
such that $AB$ forms $Q$, $AXB$ does not have an occurrence of $Q$.
Where $|A| = |\sigma_1|$ and $|B| = |\sigma_2|$ ($A$ forms $\sigma_1$ and $B$ forms $\sigma_2$).

So we have proven Lemma \ref{Lemma:pop} as well.


\end{document}